# The solution of Liouville's equation (1850, 1853) and its impact[1]

***E.M. Bogatov, S. Kichenassamy***
***Е.М. Богатов, С. Кишнассами***

Liouville's 1853 paper [2], in which he derived in closed form the general local solution of equation $u_{z\bar{z}} = e^u$, is one of the few papers from the 19th century that 21st century mathematicians routinely quote as motivation for their work. We try and understand the reasons for the enduring importance of this paper. Because it does not distinguish the real and complex domains, Liouville's paper simultaneously opened the way to a representation of the general solution of two equations, namely the elliptic equation

$$\Delta u = K e^{au} \quad \text{(L}_e\text{)}$$

and the hyperbolic equation

$$u_{xy} = K e^{au}, \quad \text{(L}_h\text{)}$$

where $K$ and $a$ are nonzero real constants. Generalizations of ($L_e$) - ($L_h$) are too numerous to be mentioned here.

Liouville announced his result in his 1850 commentary on Monge's *Applications de l'analyse à la géométrie* [3]. His original proof is lost: he only published a simpler one [2]. He obtained ($L_e$) from the expression for the Gaussian curvature of a surface in isothermal coordinates (see Lützen [4]), and obtained his general solution of ($L_h$) depending on two functions of one variable; thus, for $a = K = 1$,

$$u(x,y) = \ln \frac{2f'(x)g'(y)}{(f(x)+g(y))^2}.$$

The same formula may be interpreted as giving a real solution of ($L_e$) depending on one analytic function. The regularity and determination of the arbitrary functions was clarified rather recently (see [5, §10.6] and [6]).

Darboux [7, p.170] called ($L_h$) Liouville's equation, and Picard [8] gave the same name to ($L_e$). Picard's line of thought led to the boundary blow-up problem for ($L_e$), in relation to conformal mapping [9]-[11]. Darboux inserted ($L_h$) in a general theory of the solution of second-order PDEs in two variables. Nonetheless, Liouville's solution was not forgotten. Thus, Bianchi obtained from Liouville's solution a majorant series for the solution for what we now call the sine-Gordon equation [12, cap. XX §298, p. 214-217; cap. XXIV, §370, p. 386]. These lines of thought were apparently cut short by WWI, to be taken up again only after the birth of the theory of solitons.

Equation ($L_e$), and its counterpart in three dimensions also arose as models for isothermal gas spheres in astrophysics [13, p. 131-132], thermo-ionic emission [14, p. 50], as a special model of incompressible fluid flow as well as other applications [6, 15, 16]. Liouville's exact solution attracted considerable interest in the late 70s as an exactly solvable model of field theory [17, 18]. The contemporary emergence of the theory of solitons [19, Ch. 4] triggered a renewed interest in Bäcklund transformations, particularly the transformation relating ($L_h$) to the wave equation $v_{xy} = 0$ [20]. Liouville's solution was also applied to the motion of singularities [21, 22], and was treated by inverse scattering [23]. The modern relevance of Liouville's formula seems to be due to two facts: it is not adequately interpreted in terms of the Cauchy and Dirichlet problems; and it admits of a generalization to large classes of non-integrable equations in any number of space dimensions. On the one hand, the Cauchy problem for ($L_h$) gives a solution depending on three functions (two Cauchy data, and the equation of the curve on which they are prescribed) whereas Liouville's solution only contains two, that are essentially the two singularity data that determine a solution near blow-up [24], and this observation generalizes to non-integrable problems in any dimension; in fact, it was the starting point of the method of Fuchsian Reduction [5]. On the other hand, for ($L_e$), the Dirichlet problem requires two data (the domain and the boundary condition), whereas the boundary blow-up problem requires only one [25]-[26].

---

[1] This work is a continuation of research carried out in [1].

It is instructive to give special attention to the reception of Liouville equations in Russia and the CCCP, that seems to show at work mechanisms that may be of general significance

F. Minding, who worked since 1844 at the Dorpat University in Russia, made a great contribution to the development of the geometry of surfaces [27]. Minding was aware of all the achievements in this area but apparently, he did not know about Liouville equation, possibly because he was not interested in isothermal coordinates on surfaces. In his works, as well as in the articles of his Russian student K.M. Peterson (the founder of the Moscow geometric school), no references to the equation in either form could be found. Moscow mathematicians of that time apparently dispensed with the use of equation ($L_e$) or ($L_h$) in their research. Even now, one could argue that special coordinates detract from the identification of intrinsic quantities such as are given by tensor calculus.

By contrast, representatives of the St. Petersburg Mathematical School related to the line of P.L. Chebyshev, mentioned Liouville's results. In particular, equation ($L_e$) was present in D. Grave's article in 1927, published in the prestigious collection « *In memoriam N. I. Lobatschevskii* » in Kazan in 1927 [28, p. 26]. But it this paper the above equation was not associated with Liouville's name. It is not known how long it would have taken for this association to appear in Russian sources if the forced emigration of S.B. Bergman to Tomsk from Nazi Germany had not happened. Bergman organized a scientific seminar at the University of Tomsk on the theory of surfaces using the methods of the theory of analytic functions, which met from 1934 to 1936 [29, p. 50-51]. Apparently, it was under the influence of Bergman that B.A. Fuchs, who was one of the participants of this seminar, began to use the name "Liouville's equation" in his works using the equivalent of equation ($L_h$)

$$\frac{1}{T}\frac{\partial \log T}{\partial u \partial v} = K \qquad (L_{\log})$$

to describe a constant curvature surface's metric in 1937 [30, p. 584]. However, after 1938, Fuchs stopped referring to Liouville in connection with the use of equation ($L_{\log}$).

It so happens that in spite of the interest of Soviet mathematicians to equations of the form ($L_{\log}$) or ($L_e$), associating them with Liouville name in the mass consciousness of Soviet scientists was lost before the end of the 1970s. As a result, equation ($L_e$), obtained during the creating of the stationary theory of thermal explosion and based on the Arrhenius law, (D.I. Frank-Kamenetskii, 1939 [31, p. 740]) was also not accompanied by the references to Liouville.

It is important to note that the Dirichlet problem formed by adding to equation ($L_e$), with $K < 0$, the condition $u = 0$ on a surface $\partial\Omega$ of the form $\bar{r} = \lambda \bar{r}_0$, $(max|\bar{r}_0| = 1; \bar{r}$ is radius-vector of the point of the surface) was subjected to group-theoretic analysis by I.M. Gelfand in his famous article [32] (1959). He showed that for $K = -2$ and for symmetric domains $\Omega$, the solution to the Dirichlet problem for the equation ($L_e$) has bifurcation points when $\lambda$ reaches the threshold value $\lambda_0$, which depends on the form of $\Omega$ [32, p. 140-142].

In addition, the integral equation obtained from ($L_e$) by using Green's function for the Laplace equation turned out to be a model example of an equation of which the solutions do not belong to the space $L^p$. That is, the Liouville equation could well have been an additional stimulus for the creation of Orlicz spaces theory and its development by M.A. Krasnoselskii (for the history of this issue see [33]). Amazingly, there were no references to Liouville here either!

The situation changed only after the publication of the textbook "Modern Geometry" (1979) [34]. In Chapter 2, §13 (*Conformal form of metrics*) of this textbook, finally, equations ($L_e$) and ($L_{\log}$) were given with an indication of Liouville. At about the same time, the Liouville equation began to appear in the USSR in the works of physicists. Thus, V.A. Andreev constructed an *N*-soliton solution to equation ($L_h$), describing one moving and *N-1* stationary solitons [23]. Russian and Belarusian scientists, in relation to the two-dimensional field theory, have studied singular and regular solutions of ($L_h$), and studied generalizations of ($L_e$) and ($L_h$) etc. [21]-[22]. Soviet scientist A.M. Polyakov, known for a number of fundamental contributions to quantum field theory, found a new use for the equation ($L_e$) in high energy physics, in the quantum geometry of bosonic strings [35]: indeed, the equation of motion associated with the so-called "Liouville action"

$$S[\varphi] = C \int_\Omega \left( \frac{1}{2}\left[\left(\frac{\partial \varphi}{\partial x}\right)^2 + \left(\frac{\partial \varphi}{\partial y}\right)^2\right] + \mu^2 e^\varphi \right) dxdy; \ C = const.$$

reduces to the equation ($L_e$).

Thus, this study has shown that Liouville's paper has not been fully superseded by any of the later developments originating from it for four reasons:
1. Mathematical traditions have focused on only some aspects of it, playing an essential role in the spread (or inhibition) of mathematical ideas. The study of the Russian situation revealed the existence of several mathematical communities that received Liouville's paper differently.
2. His paper was one of the motivations of several modern theories such as nonlinear partial differential equations and functional analysis; the theory of analytical functions and conformal mapping; the theory of automorphic functions; the theory of Orlicz spaces; group-theoretical analysis; the theory of positive operators; the theory of solitons. At the same time, no single later theory can reproduce all of the known results on ($L_e$) - ($L_h$).

Furthermore, because Liouville's solution has robust features that are present in non-integrable systems in any dimension, without symmetry, Liouville's paper stimulated new theories such as Fuchsian Reduction.

**Bibliography**

1. *Богатов Е.М.* Об истории уравнения *Δu = k$e^u$* и вкладе отечественных математиков // Обозрение промышленной и прикладной математики. 2000. Том 27. Вып. 1, С. 67-69.
2. *Liouville J.* Sur l'équation aux différences partielles $\frac{d^2 \log \lambda}{du\, dv} \pm \frac{\lambda}{2a^2} = 0$ // Journal de Mathématiques Pures et Appliquées. 1853. Vol. 18. P. 71-72.
3. *Liouville J.* Sur le théorème de M. Gauss, concernant le produit des deux rayons de courbure principaux en chaque point d'une surface, Note IV to G. Monge // Applications de Analyse à la Géométrie / Paris: Bachelier. 1850. P. 583-600.
4. *Lützen J.* Joseph Liouville 1809–1882: Master of pure and applied mathematics. Springer Science & Business Media, 2012. xix+884 p.
5. *Kichenassamy S.* Fuchsian Reduction: Applications to Geometry, Cosmology and Mathematical Physics. Boston: Birkhäuser, 2007. xv+289 p.
6. *Crowdy D.G.* General solutions to the 2D Liouville equation // International Journal of Engineering Science. 1997. Vol. 35. No. 2. P. 141-149.
7. *Darboux G.* Sur les équations aux dérivées partielles du second ordre // Annales scientifiques de l'École Normale Supérieure, 1$^{re}$ série. 1870. Vol. 7. P. 163-173.
8. *Picard E.* Mémoire sur la théorie des équations aux dérivées partielles et la méthode des approximations successives // Journal de Mathématiques Pures et Appliquées. 1890. Vol. 6. P. 145-210.
9. *Poincaré H.* Les fonctions fuchsiennes et l'équation *Δu = $e^u$,* Journal de Mathématiques Pures et Appliquées. 1898. Vol. 5. No. 4. P. 137-230.

10. *Bieberbach L.* $\Delta u = e^u$ und die automorphen Funktionen (Vorläufige Mitteilung) Gesellschaft der Wissenschaften zu Göttingen.Nachrichten. Mathematisch-Physikalische Klasse 1912. S. 599-602.
11. *Bieberbach L.* $\Delta u = e^u$ und die automorphen Funktionen // Mathematische Annalen. 1916. Vol. 77. S. 173-212.
12. *Bianchi L.* Lezioni di Geometria Differenziale, vol II, 2$^{nd}$ ed., Pisa, 1903. iv+594 p.
13. *Emden R.* Gaskugeln: Anwendungen der mechanischen warmetheorie auf kosmologische und meteorologische Probleme. Leipzig und Berlin: B.G. Teubner, 1907. 497 s.
14. *Richardson O.W.* The Emission of Electricity from Hot Bodies, 2$^{nd}$ ed. London: Longmans, Green and Co., 1921. 320 p.
15. *Joseph D.D., Lundgren T.S.* Quasilinear Dirichlet problems driven by positive sources // Archive for Rational Mechanics and Analysis. 1973. Vol. 49. P. 241–269.
16. *Bandle C., Flucher M.* Harmonic Radius and Concentration of Energy; Hyperbolic Radius and Liouville's Equations $\Delta u = e^u$ and $\Delta u = u^{(n+2)/(n-2)}$// SIAM Review (Society for Industrial and Applied Mathematics). 1996. Vol. 38:2. P. 191–238.
17. *Onofri E.* An identity for *T*-ordered exponentials, with applications to Quantum Mechanics // Annals of Physics. 1976. Vol. 102. P. 371-387.
18. *Witten E.* Some exact multipseudoparticle solutions of classical Yang-Mills theory // Physical Review Letters. 1977. Vol. 38. P. 121-125.
19. *Kichenassamy S.* Nonlinear Wave Equations. New York: Dekker, 1995. xiii+288 p.
20. *Lamb G.L. Jr.* Bäcklund transformations at the turn of the century // Lecture Notes in Mathematics. 1976. Vol. 515. P. 69-79.
21. *Джорджадзе Г.П., Погребков А.К., Поливанов М. К.* Сингулярные решения уравнения □φ+(m$^2$/2)exp φ=0 и динамика особенностей // Теоретическая и математическая физика. 1979. Том 40:2. С. 221–234
22. *Pogrebkov A.K., Polivanov M.K.*, The Liouville and sinh-Gordon equations. Singular solutions, dynamics of singularities and the inverse problem method // Mathematical physics reviews. 1985. Vol. 5. P. 197–271.
23. *Андреев В.А.* Применение метода обратной задачи рассеяния к уравнению $\sigma_{xt} = e^\sigma$ // Теоретическая и математическая физика. 1976. Том 29:2, С. 213–220.
24. *Kichenassamy S., Littman W.* Blow-up Surfaces for Nonlinear Wave Equations, I & II // Communications in Partial Differential Equations. 1993. Vol. 18 (3&4). P. 431-452; P.1869-1899 .
25. *Kichenassamy S.* Boundary Blow-up and Degenerate Equations // Journal of Functional Analysis. 2004. Vol. 215 (2). P. 271-289
26. *Kichenassamy S.* Régularité du rayon hyperbolique // Comptes Rendus de l'Académie des Sciences de Paris, sér I. 2004. Vol. 338 (1). P. 13-18.
27. *Галченкова Р.И., Лумисте Ю.Г., Ожигова Е.П., Погребысский И.Б.* Фердинанд Миндинг. Л.: Наука, 1970. 225 с.
28. *Граве Д.* Плоская геометрия Эвклида как предельная для геометрии Лобачевского // In memoriam N. I. Lobatschevskii. 1927. Том 2. С. 25–36.
29. *Круликовский Н.Н.* Из истории развития математики в Томске. Томск. Изд-во Томского ун-та, 2006. 174 с.
30. *Fuchs B.* Über geodätische Mannigfaltigkeiten einer bei pseudokonformen Abbildungen invarianten Riemannschen Geometrie // Математический сборник. 1937. Том 2. №. 3. С. 567-594.
31. *Франк-Каменецкий Д.А.* Распределение температур в реакционном сосуде и стационарная теория теплового взрыва // Журнал физической химии. 1939. Том. 13. № 6. С. 738-755.
32. *Гельфанд, И. М.* Некоторые задачи теории квазилинейных уравнений // Успехи математических наук. 1959. Том. 14:2 (86). С. 87–158.

33. *Богатов Е.М.* Об истории развития нелинейных интегральных уравнений в СССР. Сильные нелинейности // Научные ведомости БелГУ. Серия Математика, Физика. 2017. Том 6. Вып. 46. С. 93-106.
34. *Дубровин Б.А., Новиков С.П., Фоменко А.Т.* Современная геометрия: методы и приложения. М.: Наука, 1979. 760 с.
35. *Polyakov A.M.* Quantum geometry of bosonic strings // Physics Letters 103B (1981), p. 207-210.